\documentclass[a4paper, 12pt]{article}
\usepackage{amssymb}
\usepackage{amsmath, amsthm, amssymb}
\usepackage{amssymb,amsmath,amsthm,url}
\usepackage{epsf,amsfonts,amsmath}

\usepackage{enumerate}
\usepackage{verbatim}

\usepackage[all]{xy}
\reversemarginpar
\def\Hom{\mbox{Hom}}
\def\mod{\mbox{mod}}
\def\Ext{\mbox{Ext}}
\def\H{\mbox{H}}
\def\Z{\mbox{Z}}
\def\DTr{\mbox{DTr}}
\def\D{\mbox{D}}
\def\Tr{\mbox{Tr}}

\theoremstyle{plain}
 \newtheorem{thm}{Theorem}[section]
 \numberwithin{equation}{section} 
 \numberwithin{figure}{section} 

 \theoremstyle{plain}
 \newtheorem{lem}[thm]{Lemma} 
 \newtheorem{prop}[thm]{Proposition} 
 \newtheorem{theorem}[thm]{Theorem}
 \newtheorem*{mythm}{Theorem}

 \theoremstyle{definition}
 \newtheorem{defn}[thm]{Definition}
 
 \newtheorem{setup}[thm]{Setup}
 
 \theoremstyle{definition}

 \theoremstyle{remark}
 \newtheorem{ex}[thm]{Example}
 
 \newtheorem{remark}[thm]{Remark}

 \newtheorem*{acknowledgement}{Acknowledgement}

\title{Existence of Auslander-Reiten sequences in subcategories}
\author{ Puiman Ng }
\date{ }

\begin{document}





\setlength{\parindent}{0pt}

\maketitle
\pagenumbering{roman}


\pagenumbering{arabic}

\begin{abstract}
This paper studies the existence of Auslander-Reiten sequences in subcategories of $\mod(\Lambda)$, where $\Lambda$ is a finite dimensional algebra over a field. The two main theorems give necessary and sufficient conditions for the existence of Auslander-Reiten sequences in subcategories.

\begin{mythm}
Let $\cal{M}$ be a subcategory of $\mod(\Lambda)$ closed under extensions and direct summands, and let $M$ be an indecomposable module in $\cal{M}$ such that $\Ext^{1}(M, \tilde{M}) \neq 0$ for some $\tilde{M}$ in $\cal{M}$. Then the following are equivalent:

\begin{itemize}
\item[(i)]
 DTr$M$ has an $\underline{\cal{M}}$-precover in the stable category $\underline{\mod}(\Lambda)$,
\item[(ii)]
 There exists an Auslander-Reiten sequence $0 \rightarrow X \rightarrow Y \rightarrow M \rightarrow 0$ in $\cal{M}$.
\end{itemize}

\end{mythm}

We also have the dual result of the above theorem. Together they strengthen the results in Auslander-Smal{\o} (\cite{AS 1}, \cite{AS 2}), and in Kleiner (\cite{Kleiner 1}).
\end{abstract}

\section{Introduction}







Let $\Lambda$ be a finite dimensional $k$-algebra over the field $k$, and let $\mod(\Lambda)$ be the category of finitely generated modules over $\Lambda$. Let $\cal{M}$ be a full subcategory of $\mod(\Lambda)$. Now let us recall a few definitions:

\begin{defn}
A morphism $g : B \rightarrow C$ in the subcategory $\cal{M}$ is said to be a right almost split morphism in $\cal{M}$ if

\begin{itemize}
\item[(i)] $g$ is not a split epimorphism,
\item[(ii)] if $h: C' \rightarrow C$ in $\cal{M}$ is not a split epimorphism, then there exists an $h': C' \rightarrow B$ such that $gh' = h$.
\end{itemize}

The notion of a left almost split morphism in the subcategory $\cal{M}$ is defined dually.
\end{defn}

\begin{defn}   \label{defnofARseq}
An exact sequence $0 \rightarrow  A \stackrel{g} \rightarrow B \stackrel{f} \rightarrow C  \rightarrow 0$ with $A, B, C$ in the subcategory $\cal{M}$ is said to be an Auslander-Reiten sequence in $\cal{M}$ if $g$ is left almost split in $\cal{M}$ and $f$ is right almost split in $\cal{M}$.
\end{defn}

Originally in \cite[Theorem 2.4]{AS 1}, Auslander and Smal{\o} developed a theory for the existence of Auslander-Reiten sequences in subcategories of $\mod(\Lambda)$. Then in \cite[Corollary 2.8]{Kleiner 1}, Kleiner gave a new proof of their existence theorem without the use of the theory of dualizing $R$-varieties. The following results strengthen the results by Auslander-Smal{\o} (\cite{AS 1}, \cite{AS 2}), and by Kleiner (\cite{Kleiner 1}), by providing necessary and sufficient conditions.

\begin{mythm}
Let $\cal{M}$ be a subcategory of $\mod(\Lambda)$ closed under extensions and direct summands, and let $M$ be an indecomposable module in $\cal{M}$ such that $\Ext^{1}(M, \tilde{M}) \neq 0$ for some $\tilde{M}$ in $\cal{M}$. Then the following are equivalent:

\begin{itemize}
\item[(i)]
 $DTrM$ has an $\underline{\cal{M}}$-precover in the injective stable category $\underline{\mod}(\Lambda)$,
\item[(ii)]
 There exists an Auslander-Reiten sequence $0 \rightarrow X \rightarrow Y \rightarrow M \rightarrow 0$ in $\cal{M}$.
\end{itemize}

\end{mythm}

\begin{mythm}

Let $\cal{L}$ be a subcategory of $\mod(\Lambda)$ closed under extensions and direct summands, and let $L$ be an indecomposable module in $\cal{L}$ such that $\Ext^{1}(\tilde{L}, L) \neq 0$ for some $\tilde{L}$ in $\cal{L}$. Then the following are equivalent:

\begin{itemize}
\item[(i)]
 $TrDL$ has an $\overline{\cal{L}}$-preenvelope in the projective stable category $\overline{\mod}(\Lambda)$,
\item[(ii)]
 There exists an Auslander-Reiten sequence $0 \rightarrow L \rightarrow B \rightarrow A \rightarrow 0$ in $\cal{L}$.
\end{itemize}

\end{mythm}

They will be proved in Theorem~\ref{leftmodulesrightARseq} and in Theorem~\ref{leftmodulesleftARseq} respectively. Their proofs are based on the theory of Auslander-Reiten triangles developed in \cite{PJ1}. The bridge between Auslander-Reiten sequences and Auslander-Reiten triangles in subcategories is shown in Lemma~\ref{ARseqARtriangle}.  \\


This paper is organized as follows: This introduction ends with some basic definitions. In Section~\ref{Setup} we describe the setup for this paper. In Section~\ref{M'precoverswitherrorterm} we introduce a weakened notion of an $\cal{M}$-precover, i.e. an $\cal{M}$-precover with error term. In Section~\ref{The stable category}, we introduce the stable categories. The main result in this section is Proposition~\ref{precoverstable}, where we show that an $\cal{M}$-precover with error term is equivalent to an $\underline{\cal{M}}$-precover in the stable category. In Section~\ref{ExistenceTheoremsforARseq}, we will prove the theorems stated above. Finally in Section~\ref{Anexample}, we will provide an example of each of the two main theorems in Section~\ref{ExistenceTheoremsforARseq}. \\

Now let $\cal T$ be some arbitrary category. Let $\cal{C}'$ be a full subcategory of $\cal T$. \\

\begin{defn}
Let $A$ and $B$ be in $\cal T$, and let $f: A \rightarrow A$ be a morphism. A morphism $\alpha: A \rightarrow B$ is \emph{right minimal} if $\alpha f = \alpha$ implies $f$ is an automorphism.
\end{defn}

\begin{defn}\label{precovercover}
A $\cal C'$-precover for an object $T$ in $\cal T$ is a morphism $c'$: $C' \rightarrow T$ for some $C'$ in $\cal C'$, such that for all $X$ in $\cal C'$, each morphism $X \rightarrow T$ factorizes through $c'$. A $\cal C'$-cover is a $\cal C'$-precover which is right minimal. The notion of a $\cal C'$-(pre)envelope is defined dually.
\end{defn}


\begin{defn}
$\cal C'$ is said to be a (pre)covering for $\cal T$ if every object in $\cal T$ has a $\cal C'$-(pre)cover.
\end{defn}


\section{Setup}  \label{Setup}
We now describe the setup for this paper. \\

Let $\Lambda$ be a finite dimensional $k$-algebra over the field $k$, and let $\mod(\Lambda)$ be the category of finitely generated $\Lambda$-left-modules. Let D be the usual duality functor $\D (-) = \Hom_{k}(-, k)$. Then $\D\Lambda$ = $\Hom_{k}(\Lambda, k)$ is the $k$-linear dual of $\Lambda$ which is a $\Lambda$-bi-module. \\

Following the notations of \cite{PJ1}, let K(Inj $\Lambda$) be the homotopy category of complexes of injective $\Lambda$-left-modules. Let $\cal T$ be the full subcategory of {K(Inj $\Lambda$)} consisting of complexes $X$ for which each $X^{i}$ is finitely generated, and where $\H^{i}(X) = 0$ for $i \gg 0$ and $\H^{i}\Hom_{\Lambda}(\D\Lambda, X) = 0$ for $i \ll 0$. $\cal T$ is triangulated. \\

Let $\cal C$ be the full subcategory of $\cal T$ which consists of injective resolutions of modules in $\mod(\Lambda)$, i.e. $\cal C$ consists of complexes in $\cal T$ of the form 

\[\xymatrix{
            \cdots  \ar[r]& 0 \ar[r] & A^{0} \ar[r] & A^{1} \ar[r] & A^{2} \ar[r] & \cdots
  }\]

where all cohomology groups other than $\H^{0}(A)$ are zero.  \\

Let $\cal{M}$ be a full subcategory of $\mod(\Lambda)$, closed under extensions and direct summands, and let $\cal C'$ consist of the injective resolutions of the $M$ in $\cal{M}$. Note that $\cal C'$ and $\cal{M}$ need not be abelian. \\


\begin{remark}\label{Cequivmodlambda}
$\cal C$ is equivalent to $\mod(\Lambda)$, by the functor $F: \mod(\Lambda) \rightarrow \cal C$ which sends $X$ in $\mod(\Lambda)$ to its injective resolution $C$ in $\cal C$. Similarly $\cal C'$ and $\cal{M}$ are equivalent. Also $\cal C'$ is (pre)covering in $\cal C$ if and only if $\cal{M}$ is (pre)covering in $\mod(\Lambda)$. \\



The setup described above could be summarized in the following diagram:

\[\xymatrix{
             \cal{M} \ar@{^{(}->}[r] \ar@{-}[d]^{\simeq} & \mod(\Lambda) \ar@{-}[d]^{\simeq} \\
             \cal C'  \ar@{^{(}->}[r] & \cal C  \ar@{^{(}->}[r] & \cal T \ar@{^{(}->}[r] &  K(\mbox{Inj} \Lambda).
  }\]

\end{remark}


Recall the notion of Auslander-Reiten triangles in subcategories (\cite[Definition 1.3]{PJ1}):

\begin{defn}\emph{(c.f. Definition~\ref{defnofARseq})} \label{defnofARtriangle}
A distinguished triangle $A \stackrel{a} \rightarrow B \stackrel{b} \rightarrow C \rightarrow$, with $A$, $B$, and $C$ in the subcategory  $\cal C'$ of $\cal T$, is defined to be an Auslander-Reiten triangle in $\cal C'$ if 
\begin{itemize}
\item[(i)] The triangle is not split,
\item[(ii)] If $A'$ is in $\cal C'$ then each morphism $A \rightarrow A'$ which is not a split monomorphism factors through $a$,
\item[(iii)] If $C'$ is in $\cal C'$ then each morphism $C' \rightarrow C$ which is not a split epimorphism factors through $b$.
\end{itemize}
\end{defn}
 
Let us now conclude the section with a little lemma. 

\begin{lem}\label{ARseqARtriangle}
$0 \rightarrow X \rightarrow Y \rightarrow Z \rightarrow 0$ is an Auslander-Reiten sequence in $\cal{M}$ if and only if $A \rightarrow B \rightarrow C \rightarrow$ is an Auslander-Reiten triangle in $\cal C'$, where $A, B, C$ are injective resolutions of $X, Y, Z$ respectively.
\end{lem}

\begin{proof}
This can be shown by standard arguments.
\end{proof}





\section{$\cal{M}$-precovers with error term}  \label{M'precoverswitherrorterm}
We are now going to introduce the notion of an $\cal{M}$-precover with error term, and discover its relationship with a $\cal C'$-precover.

\begin{setup}
Let $(-)^{\ast}$ be the functor $\Hom_{\Lambda}(-, \Lambda)$. Let $M$ be in $\cal{M}$. We will use
\[\xymatrix{
            P =  \cdots  \ar[r] &  P_{2} \ar[r] & P_{1} \ar[r] & P_{0} \ar[r] & 0 \ar[r] & \cdots
  }\]
to denote a projective resolution of $M$.

The sequence
\[\xymatrix{
            P_{1} \ar[r] & P_{0} \ar[r] & M \ar[r] & 0 
  }\]

is right exact, and the functor $(-)^{\ast}$  gives the following left exact sequence:

\[\xymatrix{
           0 \ar[r] &  M^{\ast} \ar[r] & P_{0}^{\ast} \ar[r] &  P_{1}^{\ast}.
  }\]
The \emph{transpose} Tr$M$ of $M$ is defined to be the cokernel of the map $P_{0}^{\ast} \rightarrow P_{1}^{\ast}$.
\end{setup}

\begin{remark}
Let $P_{i}$ be a projective module in $\mod(\Lambda)$. We have $\D\Lambda \otimes_{\Lambda} P_{i} \cong \Hom_{k}(\Hom_{\Lambda}(P_{i}, \Lambda), k) = \D(P_{i}^{\ast})$. Hence we have
\[\xymatrix@1@=17pt{
            \D\Lambda \otimes_{\Lambda} P = \D(P^{\ast}) =  \cdots  \ar[r] &  \D (P_{2}^{\ast}) \ar[r] & \D (P_{1}^{\ast}) \ar[r]^{d_{1}} & \D (P_{0}^{\ast}) \ar[r] & 0 \ar[r] & \cdots &.
  }\]
\end{remark}

\begin{remark}
From the exact sequence

\[\xymatrix{
           0 \ar[r] &  M^{\ast} \ar[r] & P_{0}^{\ast} \ar[r] &  P_{1}^{\ast} \ar[r] & \Tr M \ar[r] & 0,
  }\]

and the fact that the functor $\D(-)$ preserves exactness, we have the following exact sequence:

\[\xymatrix{
           0 \ar[r] &  \DTr M \ar[r] & \D (P_{1}^{\ast}) \ar[r] &  \D (P_{0}^{\ast}) \ar[r] &  \D (M^{\ast}) \ar[r] & 0.
  }\]
  
Hence DTr$M$ is the kernel of the map $d_{1}$.
\end{remark}

\begin{defn}\label{M'precoverwitherrortermdefn}
Let $M$ and $N$ be in $\cal{M}$. Then $ \nu: N \rightarrow DTrM$ is said to be an $\cal{M}$-precover with error term if for all $L$ in $\cal{M}$, each morphism $s':  L \rightarrow DTrM$ factors through $\nu$ up to an error term, i.e. there exists a morphism $\nu': L \rightarrow N$ such that $\nu \nu' - s'$ factors through $f_{2}$ in the following way: $L \rightarrow D(P_{2}^{\ast}) \stackrel{f_{2}}\rightarrow DTrM$, as indicated in the following diagram.

\[\xymatrix@1@=13pt{
& & & L \ar@<-1ex>@/_0.7pc/[dddd]_<<<<<<<<<<{{s'}} \ar@{-->}[dd]^<<<{\nu'} \ar@{.>}@/_1pc/[dddl] & &  &  && &&&&&&&& \\
 &&   &&  &&  &&   &&&&  \\
&& & N \ar'[d]^<<<{\nu}[dd]  & &  &  && &&&&&&&& \\
\Sigma^{-1}D(P^{\ast}) = \ar[rr] && D(P_{2}^{\ast}) \ar[dr]_{f_{2}}  \ar[rr]^<<<<{d_{2}} && D(P_{1}^{\ast}) \ar[rr]^{d_{1}} && D(P_{0}^{\ast}) \ar[rr] && 0 \ar[rr]  &&&& \\
&&& DTrM \ar@{^{(}->}[ur]_{j_{1}}
 }\]
 
\end{defn}

\begin{lem}\label{nullimpliesfact}
Let $M$, $N$ be in $\cal M$ and $ \nu: N \rightarrow DTrM$ be given. Let $J$ in $\cal C'$ be the injective resolution of $N$ and let $\lambda : J \rightarrow \Sigma^{-1} D(P^{\ast})$ be a chain map induced by $\nu$ as indicated in the following diagram. 

\[\xymatrix@1@=13pt{
J = \ar[dd]_{\lambda} \ar[rr] &&  0  \ar[dd] \ar[rr] && J^{0} \ar[dd]_{\lambda_{0}} \ar@{.>}@/_1pc/[ddll]_{\varphi^{0}} \ar[rr]^{a} && J^{1}\ar[dd]_{\lambda_{1}} \ar@{.>}@/_/[ddll]_{\varphi^{1}} \ar[rr] && J^{2} \ar[rr]  \ar[dd] &&&&  \\
&& & N \ar'[d]^<<{\nu}[dd] \ar@{^{(}->}[ur]_{i} & && && &&&&&&&& \\
\Sigma^{-1}D(P^{\ast}) = \ar[rr] && D(P_{2}^{\ast}) \ar[dr]_{f_{2}}  \ar[rr]^<<<<{d_{2}} && D(P_{1}^{\ast}) \ar[rr]_{d_{1}} && D(P_{0}^{\ast}) \ar[rr] && 0 \ar[rr]  &&&& \\
&&& DTrM \ar@{^{(}->}[ur]_{j_{1}}
 }\]
 
If $\lambda$ is null homotopic, then $\nu$ factorizes as $f_{2} \varphi^{0} i$ for some $\varphi^{0}$.
 
\end{lem}

\begin{proof}
Since $\lambda$ is null homotopic, $\lambda_{0} = d_{2} \varphi^{0} + \varphi^{1} a = j_{1}f_{2} \varphi^{0} + \varphi^{1} a$ for some $\varphi^{0}, \varphi^{1}$. Hence $\lambda_{0} i = j_{1}f_{2} \varphi^{0} i + \varphi^{1} a i = j_{1}f_{2} \varphi^{0} i$ since $a i = 0$. Since $\lambda_{0} i = j_{1} \nu$ (by construction), we have $ j_{1} \nu = j_{1}f_{2} \varphi^{0} i$. Since $j_{1}$ is injective, therefore $\nu = f_{2} \varphi^{0} i$.

\end{proof}

\begin{lem}\label{factimpliesnull}
Let $M$, $N$ be in $\cal M$ and $ \nu: N \rightarrow DTrM$ be given. Let $J$ in $\cal C'$ be the injective resolution of $N$ and let $\lambda : J \rightarrow \Sigma^{-1} D(P^{\ast})$ be a chain map induced by $\nu$ as indicated in the following diagram. 

\[\xymatrix@1@=13pt{
J = \ar[dd]_{\lambda} \ar[rr] &&  0   \ar[dd]  \ar[rr] && J^{0} \ar[dd]_<<<<<<<<<<{\lambda_{0}} \ar@/_1pc/[ddll]_{\varphi^{0}} \ar[rr]^{a} \ar@{->>}[dr]_<<<<{h_{1}} && J^{1}\ar[dd]^<<<<{\lambda_{1}} \ar@{.>}@/_1pc/[ddll]_<<<<<<<<{\varphi^{1}} \ar[rr]^{a^{1}} && J^{2} \ar[rr] \ar@{.>}@/_/[ddll]^<<<<<<{\varphi^{2}}  \ar[dd] &&&&  \\
&& & N \ar'[d]^<<{\nu}[dd] \ar@{^{(}->}[ur]_{i} & & \sigma N \ar@{^{(}->}[ur]_{i_{1}} \ar@{.>}@/^/[dl]^<<{g} &  && &&&&&&&& \\
\Sigma^{-1}D(P^{\ast}) = \ar[rr] && D(P_{2}^{\ast}) \ar[dr]_{f_{2}}  \ar[rr]^<<<<{d_{2}} && D(P_{1}^{\ast}) \ar[rr]_{d_{1}} && D(P_{0}^{\ast}) \ar[rr] && 0 \ar[rr]  &&&& \\
&&& DTrM \ar@{^{(}->}[ur]_{j_{1}}
 }\]
 
Suppose $\nu$ factorizes as $f_{2} \varphi^{0} i$ for some $\varphi^{0}$, then the chain map $\lambda$ is null homotopic.
 
\end{lem}

\begin{proof}
Since $\lambda_{0} i = j_{1} \nu$, we have $\lambda_{0} i = j_{1} f_{2} \varphi^{0} i = d_{2} \varphi^{0} i $. Hence $(\lambda_{0} - d_{2} \varphi^{0})i = 0$. If $\sigma N$ denotes the cokernel of $i$, then there exists a unique $g$: $\sigma N \rightarrow \D (P_{1}^{\ast})$ such that $\lambda_{0} - d_{2} \varphi^{0} = gh_{1}$. Since $\D(P_{1}^{\ast})$ is injective, there exists $\varphi^{1}: J^{1} \rightarrow \D(P_{1}^{\ast})$ such that $\varphi^{1} i_{1} = g$. Hence $\lambda_{0} - d_{2} \varphi^{0} = \varphi^{1} i_{1}h_{1}$ so $\lambda_{0} = d_{2} \varphi^{0} + \varphi^{1} i_{1}h_{1} = d_{2} \varphi^{0} + \varphi^{1} a$. Similarly, we obtain a map $\varphi^{2}: J^{2} \rightarrow \D(P_{0}^{\ast})$ such that $\lambda_{1} = d_{1} \varphi^{1} + \varphi^{2} a^{1}$.
 
\end{proof}

\begin{prop} \label{precoveringwitherrorterm} 
Let $M$ be in $\cal{M}$. Then DTr$M$ has an $\cal{M}$-precover with error term if and only if $\Sigma^{-1}$D$(P^{\ast})$ has a $\cal C'$-precover in $\cal T$. 
\end{prop}

\begin{proof}
Note that we cannot use Remark~\ref{Cequivmodlambda} since $\Sigma^{-1}$D$(P^{\ast})$ need not be in $\cal C$. We start by showing a useful diagram. 
\[\xymatrix@1@=13pt{
K = \ar@<-2ex>@/_0.7pc/[dddd]_{s} \ar@{-->}[dd]^{r}  \ar[rr] &&  0 \ar[rr] \ar[dd] && K^{0} \ar@{-->}[dd]_{r_{0}}  \ar[rr] \ar@{.>}@/_3pc/[ddddll]_<<<<<<<<<<<<<{\varphi}  && K^{1}\ar@{-->}[dd]_{r_{1}}  \ar[rr] && K^{2} \ar[rr] \ar[dd] &&&&  \\
& & & L \ar@<-1ex>@/_0.7pc/[dddd]_{s'} \ar@{-->}'[d]^<<{\nu'}[dd] \ar@{^{(}->}[ur]_{i_{l}} & &  &  && &&&&&&&& \\
J = \ar[dd]^{\lambda} \ar[rr] &&  0 \ar[rr] \ar[dd] && J^{0} \ar[dd]_{\lambda_{0}}  \ar[rr]^{a} && J^{1}\ar[dd]_{\lambda_{1}}  \ar[rr]^{a^{1}} && J^{2} \ar[rr] \ar[dd]  &&&&  \\
&& & N \ar'[d]^<<{\nu}[dd] \ar@{^{(}->}[ur]_{i} & &  &  && &&&&&&&& \\
\Sigma^{-1}D(P^{\ast}) = \ar[rr] && D(P_{2}^{\ast}) \ar[dr]_{f_{2}}  \ar[rr] && D(P_{1}^{\ast}) \ar[rr]^{d_{1}} && D(P_{0}^{\ast}) \ar[rr] && 0 \ar[rr]  &&&& \\
&&& DTrM \ar@{^{(}->}[ur]_{j_{1}}
 }\]
 
The following discussion is with reference to the diagram. \\
(\emph{only if}) Let $\nu: N \rightarrow$ DTr$M$ be an $\cal{M}$-precover with error term. Let $J$ be an injective resolution of $N$ and extend $\nu$ to a chain map $\lambda: J \rightarrow \Sigma^{-1}$D$(P^{\ast})$. We shall show that $\lambda$ is a $\cal C'$-precover. First of all, $J$ is in $\cal C'$ since $N$ is in $\cal{M}$. Suppose $K$ is in $\cal C'$ with a chain map $s: K \rightarrow \Sigma^{-1}$D$(P^{\ast})$. Then we have the induced map $s' : \Z^{0}(K) = L \rightarrow $DTr$M$. Since $L$ is in $\cal{M}$ and $\nu$ is an $\cal{M}$-precover with error term, there exists a morphism $\nu': L \rightarrow N$ such that $\nu \nu' - s' = f_{2} \varphi i_{l}$  for some $\varphi : K^{0} \rightarrow$ D$(P_{2}^{\ast})$. Extend $\nu'$ to a chain map $r: K \rightarrow J$. By Lemma~\ref{factimpliesnull}, $\lambda r - s$ is null homotopic, that is, $\lambda r = s$ in $\cal T$ and $\lambda$ is a $\cal C'$-precover. \\
(\emph{if}) Suppose $\lambda: J \rightarrow \Sigma^{-1}$D$(P^{\ast})$ is a $\cal C'$-precover. Then we get a morphism $\nu: \Z^{0}(J) = N \rightarrow$ DTr$M$. We shall show that $\nu$ is an $\cal{M}$-precover with error term. Suppose we are given $s': L \rightarrow$ DTr$M$ where $L$ is in $\cal{M}$. Extend $s'$ to a chain map $s: K \rightarrow \Sigma^{-1}$D$(P^{\ast})$ where $K$ is an injective resolution of $L$. Since $\lambda$ is a $\cal C'$-precover, there exists $r: K \rightarrow J$ such that $\lambda r = s$. This $r$ induces a homomorphism $\nu': L \rightarrow N$. By Lemma~\ref{nullimpliesfact}, $\nu \nu' - s'$ factorizes as $f_{2} \varphi^{0} i$ for some $\varphi^{0}$ and $\nu$ is therefore an $\cal{M}$-precover with error term.
\end{proof}

\section{The stable category} \label{The stable category}
In this section we study precovers in the stable category.

\begin{defn}
Let $A$ and $B$ be in $\mod(\Lambda)$. Define ${\cal I}(A, B)$ to be the set of homomorphisms from $A$ to $B$ which factor through an injective module.
\end{defn}

\begin{defn}
Let $A$ and $B$ be in $\mod(\Lambda)$. Define ${\cal P}(A, B)$ to be the set of homomorphisms from $A$ to $B$ which factor through a projective module.
\end{defn}

\begin{defn}
The (injective) stable category $\underline{\mod}(\Lambda)$ of $\mod(\Lambda)$ has the same objects as $\mod(\Lambda)$, while the morphism set $\underline{\Hom}(A,B)$ in $\underline{\mod}(\Lambda)$ is defined to be $\Hom(A,B)/{\cal I}(A, B)$ for all $A, B$ in $\underline{\mod}(\Lambda)$.
\end{defn}

\begin{defn}
The (projective) stable category $\overline{\mod}(\Lambda)$ of $\mod(\Lambda)$ has the same objects as $\mod(\Lambda)$, while the morphism set $\overline{\Hom}(A,B)$ in $\overline{\mod}(\Lambda)$ is defined to be $\Hom(A,B)/{\cal P}(A, B)$ for all $A, B$ in $\overline{\mod}(\Lambda)$.
\end{defn}

\begin{lem}\label{exactnessDP}
Let $U$ in $\mod(\Lambda)$ be a finitely generated injective $\Lambda$-module. Consider the complex $D(P^{\ast})$ from Section~\ref{M'precoverswitherrorterm}. Then
$(U,  D(P_{2}^{\ast})) \rightarrow (U,  D(P_{1}^{\ast})) \rightarrow (U,  D(P_{0}^{\ast}))$ is exact.
\end{lem}

\begin{proof}
First consider the case when $U = \D \Lambda$. Then the sequence becomes $(\D\Lambda,  \D(P_{2}^{\ast})) \rightarrow (\D\Lambda,  \D(P_{1}^{\ast})) \rightarrow (\D\Lambda,  \D (P_{0}^{\ast}))$, which is isomorphic to the sequence $(P_{2}^{\ast}, \Lambda) \rightarrow (P_{1}^{\ast}, \Lambda) \rightarrow (P_{0}^{\ast}, \Lambda)$, which is the same as the sequence $P_{2}^{\ast \ast} \rightarrow P_{1}^{\ast \ast} \rightarrow  P_{0}^{\ast \ast}$, which is isomorphic to the sequence $P_{2} \rightarrow P_{1} \rightarrow P_{0}$, which is exact. Finally, any finitely generated injective is a direct summand in a sum of copies of $\D \Lambda$. 
\end{proof}

\begin{prop}\label{precoverstable}
Let $N$ be in $\cal{M}$. Then $\nu: N \rightarrow DTrM$ is an $\cal{M}$-precover with error term in $\mod(\Lambda)$ if and only if its class $\underline{\nu} = \nu + {\cal I}(N, DTrM)$ is an $\underline{\cal{M}}$-precover in the stable category $\underline{\mod}(\Lambda)$ of $\mod(\Lambda)$.
\end{prop}




\begin{proof}

(\emph{only if}) Suppose $\nu: N \rightarrow$ DTr$M$ is an $\cal{M}$-precover with error term. We will show that $\underline{\nu} = \nu +  {\cal I} (N, \DTr M)$ is an $\underline{\cal{M}}$-precover in the stable category. Suppose we are given $\underline{s'}: L \rightarrow$ DTr$M$ in the stable category with $L$ in $\underline{\cal{M}}$, i.e. $s': L \rightarrow$ DTr$M$ in $\mod(\Lambda)$. Since $\nu$ is an $\cal{M}$-precover with error term, there exists $\nu': L \rightarrow N$ such that $\nu \nu' - s' = f_{2} \psi$ for some $\psi: L \rightarrow$ D$(P_{2}^{\ast})$. Hence $ s' = \nu \nu' - f_{2} \psi$ so $ \underline{s'}$ = $\underline{\nu}$ $\underline{\nu'}$ - $\underline{f_{2}}$ $\underline{\psi}$ = $\underline{\nu}$ $\underline{\nu'}$, since $\underline{f_{2}}$ = $\underline{0}$.

(\emph{if}) Suppose $\underline{\nu}$ is an $\underline{\cal{M}}$-precover in the stable category. We will show  $\nu$ is an $\cal{M}$-precover with error term. Suppose we are given $s': L \rightarrow$ DTr$M$ with $L$ in $\cal{M}$. Consider its class $ \underline{s'}$ in the stable category. Since $\underline{\nu}$ is an $\underline{\cal{M}}$-precover in the stable category, we have $\underline{\nu'}: L \rightarrow N$ such that $ \underline{s'}$ =  $\underline{\nu}$ $\underline{\nu'}$ for some $\underline{\nu'}$, i.e. $\nu \nu' - s'$ factors through an injective $U$, say $\nu \nu' - s' = u_{2} u_{1} $. We would like $\nu \nu' - s'$, however, to factor through D$(P_{2}^{\ast}) \stackrel{f_{2}} \rightarrow$ DTr$M$ instead. 
\[\xymatrix{
           L \ar[rr]^{\nu \nu' - s'} \ar[dr]^{u_{1}} \ar@{-->}[dddr] &   & DTrM \\
           &  U \ar[ur]^{u_{2}} \ar@{-->}[dd] & \\
           & & \\
           &  D(P_{2}^{\ast})       \ar[uuur]_{f_{2}} &
             }\]
             
Consider the morphism $j_{1}u_{2}: U \rightarrow$D$(P_{1}^{\ast})$, where $j_{1}$ is as in Definition~\ref{M'precoverwitherrortermdefn}.

\[\xymatrix{ & U \ar[d]^{j_{1}u_{2}} \ar@{-->}[dl]_{g} \ar[dr]^{0} & \\
D(P_{2}^{\ast}) \ar[r]_{d_{2}}  & D(P_{1}^{\ast}) \ar[r]_{d_{1}}  & D(P_{0}^{\ast}) 
    }\]
    
Since $d_{1}j_{1}u_{2} = 0$, by Lemma~\ref{exactnessDP}, there exists $g$ such that $d_{2} g = j_{1}u_{2}$, which gives $j_{1} f_{2} g = j_{1}u_{2}$. Since $j_{1}$ is injective, therefore $f_{2} g = u_{2}$ and $\nu \nu' - s = u_{2} u_{1} = f_{2} g u_{1}$. 
\end{proof}

\section{Existence of Auslander-Reiten sequences in subcategories} \label{ExistenceTheoremsforARseq}

Let me restate Theorem 3.1 of \cite{PJ1} here: Let $C$ be in $\cal C$ and let $X \rightarrow Y \rightarrow C \rightarrow$ be an Auslander-Reiten triangle in $\cal T$. Then the following are equivalent:

\begin{itemize}
\item[(i)] $X$ has a $\cal C$-cover of the form $A \stackrel{\alpha} \rightarrow X$.
\item[(ii)] There is an Auslander-Reiten triangle $A \rightarrow B \rightarrow C \rightarrow$ in $\cal C$.
\end{itemize}



Now we state the existence theorem for (right) Auslander-Reiten sequences in subcategories: \\ 

\begin{theorem}\label{leftmodulesrightARseq}
Let $\cal{M}$ be a subcategory of $\mod(\Lambda)$ closed under extensions and direct summands, and let $M$ be an indecomposable module in $\cal{M}$ such that $\Ext^{1}(M, \tilde{M}) \neq 0$ for some $\tilde{M}$ in $\cal{M}$. Then the following are equivalent:

\begin{itemize}
\item[(i)]
 DTr$M$ has an $\underline{\cal{M}}$-precover in the injective stable category $\underline{\mod}(\Lambda)$,
\item[(ii)]
 There exists an Auslander-Reiten sequence $0 \rightarrow X \rightarrow Y \rightarrow M \rightarrow 0$ in $\cal{M}$.
\end{itemize}

\end{theorem}

\begin{proof}
(i) $\Rightarrow$ (ii):
Let $P$ and $C$ be a projective and an injective resolution of $M$. By Proposition~\ref{precoverstable}, DTr$M$ has an $\cal{M}$-precover with error term, and then by Proposition~\ref{precoveringwitherrorterm}, $\Sigma^{-1} \D (P^{\ast})$ has a $\cal C'$-precover. Using Theorem 4.6 in \cite{PJ1}, there exists an Auslander-Reiten triangle $A \rightarrow B \rightarrow C \rightarrow$ in $\cal C'$, and finally an Auslander-Reiten sequence  $0 \rightarrow \H^{0}(A) \rightarrow \H^{0}(B) \rightarrow  \H^{0}(C) \rightarrow 0$ in $\cal{M}$, where $M$ is retrieved through the isomorphism $\H^{0}(C) \cong M$ (Lemma~\ref{ARseqARtriangle}). \\

(ii) $\Rightarrow$ (i):
Let $P$ and $C$ be a projective and an injective resolution of $M$. Following the argument in Theorem 4.6 in \cite{PJ1}, there exists an Auslander-Reiten triangle $\Sigma^{-1} \D (P^{\ast}) \rightarrow Y \rightarrow  C \rightarrow$ in $\cal T$. Since there exists an Auslander-Reiten sequence $0 \rightarrow X \rightarrow Y \rightarrow M \rightarrow 0$ in $\cal{M}$, therefore by Remark~\ref{ARseqARtriangle}, there exists an Auslander-Reiten triangle $A \rightarrow B \rightarrow C \rightarrow $ in $\cal C'$. By Theorem 3.1 in \cite{PJ1}, $\Sigma^{-1} \D (P^{\ast})$ has a $\cal C'$-precover. By Proposition~\ref{precoveringwitherrorterm}, DTr$M$ has an $\cal{M}$-precover with error term. Finally by Proposition~\ref{precoverstable}, DTr$M$ has an $\underline{\cal{M}}$-precover in the stable category $\underline{\mod}(\Lambda)$.
\end{proof}

Before we give the existence theorem for (left) Auslander-Reiten sequences in subcategories, we need the following:

\begin{lem}\label{ARseqdual}
Let $\cal{L}$ be a subcategory of $\mod(\Lambda)$. Then the following are equivalent:

\begin{itemize}
\item[(i)]
$0 \rightarrow X \rightarrow Y \rightarrow Z \rightarrow 0$ is an Auslander-Reiten sequence in $\cal{L}$, 
\item[(ii)]
$0 \rightarrow DZ \rightarrow DY \rightarrow DX \rightarrow 0$ is an Auslander-Reiten sequence in $D\cal{L}$. 
\end{itemize}

\end{lem}

\begin{proof}
Refer to the remark after \cite[Proposition V.1.13]{ARS 1}.
\end{proof}

\begin{prop}\label{precoverpreenvelopeinstablecategory}
Let $\cal{L}$ be a subcategory of $\mod(\Lambda)$. Denote $D\cal{L}$ by $\cal{M}$. Let $L$ be in $\cal{L}$.
Then the following are equivalent:

\begin{itemize}
\item[(i)]
$DL$ has an $\underline{\cal{M}}$-precover in the injective stable category $\underline{\mod}(\Lambda^{op})$, 
\item[(ii)]
$L$ has an $\overline{\cal{L}}$-preenvelope in the projective stable category $\overline{\mod}(\Lambda)$.
\end{itemize}
\end{prop}

\begin{proof}





This can be shown by standard arguments.
\end{proof}





\begin{prop}\label{dualversionofleftmodulesrightARseq}
Let $\cal{L}$ be a subcategory of $\mod(\Lambda)$ closed under extensions and direct summands. Denote $D\cal{L}$ by $\cal{M}$. Let $L$ be an indecomposable module in $\cal{L}$ such that $\Ext^{1}(DL, D\tilde{L}) \neq 0$ for some $\tilde{L}$ in $\cal{L}$. Then the following are equivalent:

\begin{itemize}
\item[(i)]
 $DTrDL$ has an $\underline{\cal{M}}$-precover in the injective stable category $\underline{\mod}(\Lambda^{op})$,
\item[(ii)]
 There exists an Auslander-Reiten sequence $0 \rightarrow DA \rightarrow DB \rightarrow DL \rightarrow 0$ in $\cal{M}$, where $A$ and $B$ are in $\cal{L}$.
\end{itemize}
\end{prop}

\begin{proof}
Since $\cal{M}$ = $\D\cal{L}$, and $\cal{L}$ is closed under extensions and direct summands, therefore $\cal{M}$ is a subcategory of $\mod(\Lambda^{op})$ closed under extensions and direct summands. Since $L$ is an indecomposable module in $\cal{L}$, $\D L$ is an indecomposable module in $\cal{M}$. The rest follows from the right module version of Theorem~\ref{leftmodulesrightARseq}. 
\end{proof}

Finally we have the dual of Theorem~\ref{leftmodulesrightARseq}:
\begin{theorem}\label{leftmodulesleftARseq}
Let $\cal{L}$ be a subcategory of $\mod(\Lambda)$ closed under extensions and direct summands, and let $L$ be an indecomposable module in $\cal{L}$ such that $\Ext^{1}(\tilde{L}, L) \neq 0$ for some $\tilde{L}$ in $\cal{L}$. Then the following are equivalent:

\begin{itemize}
\item[(i)]
 $TrDL$ has an $\overline{\cal{L}}$-preenvelope in the projective stable category $\overline{\mod}(\Lambda)$,
\item[(ii)]
 There exists an Auslander-Reiten sequence $0 \rightarrow L \rightarrow B \rightarrow A \rightarrow 0$ in $\cal{L}$.
\end{itemize}

\end{theorem}

\begin{proof}
Let $\cal{M}$ = $\D\cal{L}$. By Proposition~\ref{precoverpreenvelopeinstablecategory}, TrD$L$ has an $\overline{\cal{L}}$-preenvelope in $\overline{\mod}(\Lambda )$ if and only if DTrD$L$ has an $\underline{\cal{M}}$-precover in $\underline{\mod}(\Lambda^{op} )$. Note that $\Ext^{1}(\tilde{L}, L) \neq 0$ if and only if $\Ext^{1}(\D L, \D\tilde{L}) \neq 0$. Hence the result follows from Proposition~\ref{dualversionofleftmodulesrightARseq}, with the help of Lemma~\ref{ARseqdual}.
\end{proof}

\section{An example} \label{Anexample}
In this section we are going to give an example of Theorem~\ref{leftmodulesrightARseq} and of Theorem~\ref{leftmodulesleftARseq}. 

\begin{ex}
Let $\Lambda$ be a representation-infinite hereditary algebra. Let $\cal{M}$ be any full subcategory of $\mod(\Lambda)$ which consists of postprojective modules and is closed under extensions and direct summands, see \cite[Definition VIII.2.2]{Assem 1}. Let $M$ be an indecomposable module in $\cal{M}$ such that $\Ext^{1}(M, \tilde{M}) \neq 0$ for some $\tilde{M}$ in $\cal{M}$. Then there exists an Auslander-Reiten sequence $0 \rightarrow X \rightarrow Y \rightarrow M \rightarrow 0$ in $\cal{M}$.

\end{ex}

\begin{proof}
By \cite[Lemma VIII.2.5]{Assem 1}, we know that there are only finitely many indecomposable modules in $\cal{M}$ which have non zero maps to DTr$M$. Hence DTr$M$ has an $\cal{M}$-precover in $\mod(\Lambda)$ and therefore an $\underline{\cal{M}}$-precover in the stable category $\underline{\mod}(\Lambda)$. The existence of the Auslander-Reiten sequence $0 \rightarrow X \rightarrow Y \rightarrow M \rightarrow 0$ in $\cal{M}$ follows from Theorem~\ref{leftmodulesrightARseq}. 
\end{proof}

\begin{remark}
Dually, let $\cal{L}$ be the full subcategory of $\mod(\Lambda)$ consisting of preinjective modules over $\Lambda$, and let $L$ be an indecomposable module in $\cal{L}$ such that $\Ext^{1}(\tilde{L}, L) \neq 0$ for some $\tilde{L}$ in $\cal{L}$. 
Then the existence of the Auslander-Reiten sequence $0 \rightarrow L \rightarrow Y \rightarrow X \rightarrow 0$ in $\cal{L}$ follows from Theorem~\ref{leftmodulesleftARseq}.
\end{remark}

\medskip

\begin{acknowledgement}
 The author would like to express her thankfulness for the financial assitance from the School of Mathematics and Statistics, Newcastle University, for the Newcastle University International Postgraduate Scholarship (NUIPS), for the Overseas Research Students Awards Scheme (ORSAS) Award, and for the Croucher Foundation Scholarship, which have enabled her study in Newcastle University. The author would also like to express her special gratitude to her supervisor, Peter J\o rgensen, for whose ideas he generously imparts, that have come into their right signification through his patient guidance.
\end{acknowledgement}










\begin{thebibliography}{1}

\bibitem{Assem 1}
Ibrahim Assem, Andrzej Skowronski, and Daniel Simson,
\newblock Elements of the Representation Theory of Associative Algebras: Volume 1: Techniques of Representation Theory,
\newblock London Mathematical Society Student Texts 65, Paperback, Cambridge University Press, Cambridge, 2006.

\bibitem{ARS 1}
M. Auslander,  Idun Reiten, and S. O. Smal\o ,
\newblock Representation Theory of Artin Algebras,
\newblock Cambridge Studies in Advanced Mathematics 36, Cambridge University Press, Cambridge, 1995.

\bibitem{AS 1}
M. Auslander and S. O. Smal\o ,
\newblock Almost Split Sequences in Subcategories,
\newblock J. Algebra 69(2), 426 - 454 (1981).

\bibitem{AS 2}
M. Auslander and S. O. Smal\o ,
\newblock Addendum to ``Almost Split Sequences in Subcategories",
\newblock J. Algebra 71 (1981), 592 - 594.


\bibitem{PJ1} 
Peter J\o rgensen,
\newblock Auslander-Reiten Triangles in Subcategories,
\newblock J. K-theory 3 (2009), 583-601.


\bibitem{Kleiner 1}
Mark Kleiner,
\newblock Approximations and Almost Split Sequences in Homologically Finite Subcategories,
\newblock J. Algebra 198, 135 - 163 (1997).

\end{thebibliography}

\medskip

School of Mathematics and Statistics,
Newcastle University,
Newcastle upon Tyne, NE1 7RU,
United Kingdom

\medskip

e-mail address: p.i.ng@newcastle.ac.uk

\end{document}